\documentclass[11pt,a4paper]{article}

\usepackage{amsmath,amssymb}
\newcommand{\R}{\mathbb{R}}

\newcommand{\p}{\partial}

\newtheorem{lem}{Lemma}

\newcommand{\e}{\epsilon}
\newcommand{\va}{\varphi}

\newtheorem{theorem}{Theorem}

\newtheorem{proposition}{Proposition}

\newtheorem{remarka}{Remark}

\title{Existence of global strong solution for the compressible Navier-Stokes equations with degenerate viscosity coefficients in 1D}
\author{Boris Haspot  \thanks{Universit\'e Paris Dauphine, PSL Research University, Ceremade, Umr Cnrs 7534, Place du Mar\' echal De Lattre De Tassigny 75775 Paris cedex 16 (France), haspot@ceremade.dauphine.fr } \footnote{ the author certify that the general content of the manuscript, in whole or in part, is not submitted, accepted, or published elsewhere, including conference proceedings. The author also certify that the general content of the manuscript, in whole or in part, will not be submitted, accepted, or published elsewhere, including conference proceedings, while it is under consideration or in production.}}
\date{}
\begin{document}
\maketitle


\begin{abstract}
We consider Navier-Stokes equations for compressible viscous fluids in the one-dimensional case. We prove the existence of global strong solution with large initial data for compressible Navier Stokes equation with viscosity coefficients of the form $\p_x(\rho^\alpha\p_x u)$ with $\frac{1}{2}<\alpha\leq 1$ (it includes in particular the important physical case of the viscous shallow water system when $\alpha=1$). The key ingredient of the proof relies to a new formulation of the compressible equations involving a new effective velocity $v$ (see \cite{cras,para,CPAM,CPAM1}) such that the density verifies a parabolic equation. We estimate $v$ in $L^\infty_{t,x}$ norm which enables us to control the $L^\infty_{t,x}$ norm of $\frac{1}{\rho}$ by using the maximum principle.
\end{abstract}
\section{Introduction}
In this paper we wish to investigate the existence of global strong solutions of the following Navier-Stokes equations for compressible isentropic flow:
\begin{equation}
\begin{cases}
\begin{aligned}
&\p_t\rho+\p_x(\rho u)=0,\\
&\p_t(\rho u)+\p_x(\rho u^2)-\p_x(\mu (\rho)\p_x u)+\p_x P(\rho)=0.
\end{aligned}
\end{cases}
\label{1}
\end{equation}
with possibly degenerate viscosity coefficient $\mu(\rho)\geq 0$.\\
Throughout the paper, we will assume that the pressure $P(\rho)$ verifies a $\gamma$ type law:
\begin{equation}
P(\rho)=\rho^\gamma,\;\gamma>1.
\end{equation}
Following the idea of \cite{BrDeZa,cras,para,CPAM,CPAM1}, setting $v=u+\p_x \va(\rho)$ with $\va'(\rho)=\frac{\mu(\rho)}{\rho^2}$ we can rewrite the system (\ref{1}) as follows \footnote{We refer to the Appendix for more details on the computations.}:
\begin{equation}
\begin{cases}
\begin{aligned}
&\p_t\rho-\p_{x}(\frac{\mu(\rho)}{\rho}\p_x\rho)+\p_x(\rho v)=0,\\
&\rho \p_t v+\rho u\p_x v+\p_x P(\rho)=0.
\end{aligned}
\end{cases}
\label{11}
\end{equation}
The important point to note here is that this change of unknown is true for any viscosity coefficients (including in particular the important case of constant viscosity coefficients). There is no limitation as it is the case in dimension $N\geq 2$, indeed we recall that in this case the BD entropy involves an algebraic relation between the coefficients $\mu$ and $\lambda$ which precludes in particular the case of constant viscosity coefficients (see \cite{BrDeZa,para} for more details).\\
We emphasize in addition that the introduction of the effective velocity $v$ allows to transform the system (\ref{1}) into a parabolic equation on the density and a transport equation on the velocity (we will see in the sequel that $v$ is not so far to verify a damped transport equation). It seems surprising to observe that contrary to $u$ which has a parabolic behavior, $v$  has a hyperbolic behavior. Roughly speaking, the compressible Navier Stokes equations in the one dimensional-case can be seen as the compressible Euler equations with a viscous regularizing term on the density of the type $-\p_{x}(\frac{\mu(\rho)}{\rho}\p_x\rho)$. Note that it is less obvious in dimension $N\geq 2$ since the momentum equation of (\ref{11}) has in addition a term of the form ${\rm div}(\mu(\rho){\rm curl}v)$, see \cite{para}).\\
A large amount of literature is dedicated to the study of the compressible Navier-Stokes equations with
 constant viscosity case, however physically the viscosity of a gas depends on the temperature and  on the density (in the isentropic case). Let us mention the case of the Chapman-Enskog viscosity law (see \cite{CC70}) or the case of monoatomic gas ($\gamma=\frac{5}{3}$) where $\mu(\rho)=\rho^{\frac{1}{3}}$. More generally, the viscosity coefficient $\mu(\rho)$ is expected to vanish as a power of the density $\rho$ on the vacuum. In this paper we are going to deal with degenerate viscosity coefficients which can be written under the form $\mu(\rho)=\mu\rho^\alpha$ with $\frac{1}{2}<\alpha\leq 1$. It is worth pointing out that the case $\alpha=1$ corresponds to the so called viscous shallow water system. This system with friction has been derived by Gerbeau and Perthame in \cite{GP} from the Navier-Stokes system with a free moving boundary in the shallow water regime at the first order (it corresponds to a small shallowness parameter). This derivation relies on the hydrostatic approximation where the authors follow the role of viscosity and friction on the bottom. 
 \\
 \\
We are going now to recall some results on the existence of solutions for the one-dimensional case when the viscosity coefficient is constant positive. The existence of global weak solutions was first obtained by Kazhikhov and Shelukin \cite{KS77} for smooth enough data close to the equilibrium (in particular the initial density $\rho_0$ is bounded away from zero). The case of discontinuous data (still bounded away from zero) was studied by Shelukin \cite{She82,She83,She84} and then by Serre \cite{Ser86a,Ser86b} and Hoff \cite{Hof87}. First results dealing with vanishing initial density were also obtained by Shelukin  \cite{She86}. In \cite{Hof98}, Hoff extends the previous results by proving the existence of global weak solution with large discontinuous initial data having different limits at $x=\pm \infty$. In passing let us mention that the existence of global weak solution in any dimension $N\geq 2$ has been proved for the first time by Lions in \cite{Lio98} and the result has been later refined by Feireisl et al (\cite{FNP01} and \cite{Fei04}, see also \cite{BJ}). Concerning the uniqueness of the solution, Solonnikov  in \cite{Sol76} obtained the existence of strong solution for smooth initial data in finite time. However, the regularity may blow up when the density approaches from the vacuum. A natural question is to understand if the vacuum may  occur in finite time. Hoff and Smoller (\cite{HS01}) showed that for any weak solution of the Navier-Stokes equation 
in the one-dimensional case, the density remains strictly positive all along the time provided that no vacuum states exist initially. The existence of global strong solution with large initial data for initial density far away from the vacuum has been proved for the first time by Kanel \cite{Ka}  (see also \cite{Hof87}).
\\
\\
Concerning the case with degenerate viscosity coefficients,  it seems more delicate to estimate all along the time the $L^\infty$ norm of $\frac{1}{\rho}$ when we assume that initially $\frac{1}{\rho_0}$ belongs to $L^\infty$. In other words, it is possible that vacuum arises in finite time. This is in particular the main obstruction in order to prove the existence of global strong solution for degenerate viscosity coefficients. This question is at present far from being solved. The main theorem of this paper answers to the case $\mu(\rho)=\mu\rho^\alpha$ with $\frac{1}{2}<\alpha\leq 1$.
We emphasize that the 1D compressible Navier-Stokes equations with degenerate viscosity 
has also been considered in the opposite case where we assume that the initial density is compactly supported (see \cite{LXY98},  \cite{Mu,Mu1}, \cite{OMNM02}, \cite{YYZ01}, \cite{YZ02}, \cite{JXZ} and \cite{LLX}). The authors are concerned by the evolution of the free boundary delimiting the vacuum. \\
As explained previously, in our case the crucial fact is that $\frac{1}{\rho}$ remains bounded in $L^\infty$ norm all along the time in order to take into account the parabolic effects on the velocity $u$. This point is essential if we wish to prove global uniqueness results. The first results are due to Mellet and Vasseur who proved in \cite{MV} the existence of global strong solution with large initial data when $\mu(\rho)$ verifies the following condition:
 $$
 \begin{cases}
 \begin{aligned}
 &\mu(\rho)\geq\nu\rho^{\alpha}\;\;\forall \rho\leq 1\;\;\mbox{for some}\;\alpha\in[0,\frac{1}{2}),\\
 &\mu(\rho)\geq\nu\;\;\forall \rho\geq 1.
 \end{aligned}
 \end{cases}
 $$
The key tool of the proof consists in controlling the $L^\infty$ norm of $\frac{1}{\rho}$ by using an entropy derived by Bresch and  Desjardins in \cite{BD} (see also the interesting paper \cite{BrDeZa}) in the multi-dimensional case ( in the one-dimensional case, a similar inequality was introduced earlier by Kanel in \cite{Ka} and Vaigant \cite{Vai90} for flows with constant viscosity, see also  Shelukin \cite{She98}). Indeed it allows to Mellet and Vasseur to control in $L^\infty(L^2)$ norm the quantity $\p_x\rho^{\alpha-\frac{1}{2}}$, from Sobolev embedding they deduce that  $\frac{1}{\rho}$ belongs to $L^\infty_T(L^\infty(\R))$ for any $T>0$. Similar arguments allow to control the $L^\infty$ norm of the density $\rho$. It is then sufficient for proving that the solution of V. Solonnikov in \cite{Sol76} can be extended for any time $T$ (or in other word that there is no blow up in finite time). We would like also to mention that Mellet and  Vasseur proved in \cite{MV06} the stability of the global weak solution for compressible Navier-Stokes equation with viscosity coefficient verifying the so called BD entropy (see \cite{BDa,BD}) in dimension $N= 2,3$ (we refer also to \cite{BD,BrDeZa} when we add friction terms). Let us mention in particular that the case $\mu(\rho)=\mu\rho$ with $\mu>0$ and $\lambda(\rho)=0$ verifies the algebraic relation related to the BD entropy discovered in \cite{BD}, it corresponds here to the so called viscous shallow water system.
For $N=2,3$ the important problem of the existence of global weak solutions has been recently resolved independently by Vasseur and Yu \cite{V,V1} and Li and Xin in \cite{Li}.
\\
The goal of this paper consists in extending the result of \cite{MV} to the case where $\mu(\rho)=\mu\rho^\alpha$ with $\frac{1}{2}<\alpha\leq 1$ and $\mu>0$. To do this, we are going to study the system (\ref{11}) and to prove that $v$ remains in $L^\infty_T(L^\infty)$ for any $T>0$. Indeed we can observe that the effective velocity $v$ verifies a damped transport equation with a remainder term bounded in $L^2_T(L^\infty)$. We can then use the maximum principle on the first equation of (\ref{11}) in order to control $\frac{1}{\rho}$ in $L^\infty$ norm. In particular we deduce that the solution of Solonnikov in \cite{Sol76} does not blow up for any time $T>0$ which is sufficient to show the existence of global strong solution.\\
In the next section we state our main result. Section \ref{section2} deals with the proof of  the theorem \ref{theo1} and we postpone an Appendix in order to explain the equivalence between the system (\ref{1}) and the system (\ref{11}).
\section{Main result}
Following Hoff in \cite{Hof98}, we work with positive initial data having positive limits at $x=\pm\infty$. We fix constant positive density $\rho_{+}>0$ and $\rho_{-}>0$ and a smooth monotone function $\bar{\rho}(x)$ such that:
\begin{equation}
\bar{\rho}(x)=\rho_{\pm}\;\;\mbox{when}\;\;x_{\pm}\geq 1,\;\;\bar{\rho}(x)\geq c>0\;\forall x\in\R.
\end{equation}
We assume in addition that $\p_x\bar{\rho}$ belongs to $L^2(\R)\cap L^\infty(\R)$. In the sequel in order to simplify the proof, we will deal with viscosity coefficients of the form $\mu(\rho)=\mu\rho^{\alpha}$ with $\alpha>0$. In addition we assume that there exists $C>0$ such that:
\begin{equation}
\mu(\rho)\leq C+CP(\rho)\;\;\forall\rho\geq 0. 
\label{2.4}
\end{equation}
This condition will be verified in the assumptions of the Theorem \ref{theo1}. Our main theorem states as follows.\begin{theorem}
Assume that $\mu(\rho)=\rho^\alpha$ with $\frac{1}{2}<\alpha\leq 1$ and $0<\e<\frac{1}{4}$.
Assume that  $P(\rho)=a\rho^{\gamma}$ with $\gamma\geq \alpha+\frac{1}{2}+\e$. The initial data $\rho_0$ and $u_0$ satisfy:
\begin{equation}
\begin{aligned}
&0<\alpha_0\leq\rho_0(x)\leq\beta_0<+\infty,\\
&\rho_0-\bar{\rho}\in H^1(\R),\\
&u_0 \in H^1(\R),\\
&\p_x\rho_0\in L^\infty,
\end{aligned}
\label{2.5}
\end{equation}
for some constant $\alpha_0$ and $\beta_0$. Then there exists a global strong solution $(\rho,u)$ of system (\ref{1}) on $\R^+\times\R$ such that for every $T>0$:
\begin{equation}
\begin{aligned}
&\rho-\bar{\rho}\in L^\infty(0,T;H^1(\R)),\\
&u\in L^\infty(0,T;H^1(\R))\cap L^2(0,T,H^2(\R)),\\
&v\in L^\infty_T(L^\infty(\R)).
\end{aligned}
\label{2.55}
\end{equation}
Moreover for every $T>0$, there exists constant $\alpha_1(T)$ and $\beta(T)$ such that
\begin{equation}
0<\alpha_1(T)\leq\rho(t,x)\leq\beta(T)<+\infty\;\;\forall (t,x)\in(0,T)\times\R.
\label{2.555}
\end{equation}
\label{theo1}
\end{theorem}
\begin{remarka}
This result extends the work \cite{MV} to the classical viscous shallow water system and more generally to the case $\frac{1}{2}<\alpha\leq 1$.\\
We emphasize that it would be possible to consider initial velocity $u_0$ verifying $(u_0-\bar{u})\in H^1(\R)$ with $\bar{u}$ a regular function having different limits at the infinity.\\
It should be also possible to weaken the regularity assumptions on the initial data by working in critical space for the scaling of the equations (see \cite{JDE,JDE2}).
\end{remarka}
\begin{remarka}
We emphasize that our conditions on $\gamma$ allow to deal with the so-called viscous Saint Venant equation derived by Gerbeau and Perthame in \cite{GP} from the Navier-Stokes equations with free boundary where $P(\rho)=\rho^2$, $\gamma=2$.\\
We would like to precise why the theorem has restriction on the range of $\alpha$ and $\gamma$. First the fact that $\alpha>\frac{1}{2}$ allows to show that the density remains bounded in $L^\infty$ norm all along the time. Indeed from Sobolev embedding this is essentially  a consequence of the fact that $\p_x(\rho^{\alpha-\frac{1}{2}})$ belongs to $L^\infty_T(L^2)$ for any $T>0$. Now if we wish to estimate the $L^\infty$ norm of $\frac{1}{\rho}$, it seems natural to use the maximum principle on the mass equation of (\ref{11}). To do this we need to verify that the diffusion term $\frac{\mu(\rho)}{\rho}=\rho^{\alpha-1}$ is bounded away from zero and this requires that $\alpha\leq 1$ since $\rho$ belongs to $L^\infty_{t,x}$.\\
In addition in order to apply the maximum principle, we must be able to estimates $v\in L^\infty(L^p)$ for $p$ sufficiently large.
Since $v$ verify a damped transport equation with a remainder term of the form $\frac{P(\rho)}{\mu(\rho)}u$, it remains essentially to control $\frac{P(\rho)}{\mu(\rho)}u$ in $L^1_T(L^\infty)$ for any $T>0$. We will show that $\rho^{\frac{1}{2}+\e}u$ remains bounded in $L^2_T(L^\infty)$ for any $T>0$. Now we can observe that $\frac{P(\rho)}{\mu(\rho)}u=\rho^{\gamma-\alpha-\frac{1}{2}-\e}\rho^{\frac{1}{2}+\e}u$, it explains why we assume that $\gamma-\alpha-\frac{1}{2}-\e\geq 0$ because $\rho^{\gamma-\alpha-\frac{1}{2}-\e}$ can be estimated via the $L^\infty$ norm of the density $\rho$.\\
We think that our assumption $\gamma\geq\alpha+\frac{1}{2}+\e$ is not optimal, in the proposition \ref{proposition4.2.2b}. It seems that the choice $\beta=\frac{\alpha}{2}+\e$ is better (in this case the condition on $\gamma$ would be $\gamma\geq\frac{3\alpha}{2}+\e$). In order to do this we must obtain estimates on $\rho^{\frac{1}{p}}u$ in $L^\infty_T(L^p)$ with $1\leq p<2$ for any $T>0$. A good option to prove this is maybe to deal with weight on the velocity $u_0$. 
\\
We would like to mention that in the result of Mellet and Vasseur \cite{MV}, the coefficient $\gamma$ verifies the less restrictive condition $\gamma>1$. This is due to the fact that the control on $\frac{1}{\rho}$ in $L^\infty$ direct is more direct and is a consequence of the entropy on $v$.\\
\end{remarka}
When the viscosity coefficient $\mu(\rho)$ satisfies:
\begin{equation}
\mu(\rho)\geq\nu>0\;\;\forall\rho\geq 0,
\label{2.6}
\end{equation}
the existence of strong solution with large initial data in finite time is classical (see \cite{Sol76}). 
\begin{proposition} Let $(\rho_0, u_0)$ satisfy (\ref{2.5}) and assume that $\mu$ satisfies (\ref{2.6}) there exists $T_0$ depending on $\alpha_0,\beta_0,\|\rho_0-\bar{\rho}\|_{H^1}$ and $\|u_0\|_{H^1}$ such that (\ref{1}) has a unique solution $(\rho,u)$ on $(0,T_0)$ satisfying:
$$
\begin{aligned}
&\rho-\bar{\rho},u \in L^\infty(0,T_1,H^1(\R)),\;\p_t\rho\in L^2((0,T_1)\times\R),\\
&u\in L^2(0,T_1,H^2(\R)),\;\p_t u\in  L^2((0,T_1)\times\R)
\end{aligned}
$$
for all $T_1<T_0$.\\
Moreover, there exist some $\alpha(T)>0$ and $\beta(T)<+\infty$ such that $\alpha(t)\leq \rho(,x)\leq\beta(t)$ for all $t\in (0,T_0)$.
\label{Solo}
\end{proposition}
\begin{remarka}
It is important to mention that if $\rho_0-\bar{\rho}$ and $u_0$ are in $H^s(\R)$ with $s\geq 1$ then the regularity $H^s$ is preserved on $(0,T_0)$ and we have $\rho-\bar{\rho},\,u \in L^\infty(0,T,H^s(\R))$ for any $T<T_0$.
\end{remarka}
\begin{remarka}
We mention that it could be possible to extend this theorem to the case of critical initial data for the scaling of the equations (see \cite{M3AS,JDE,JDE2} in dimension $N\geq 2$). In particular it would allow to improve the theorem \ref{theo1} in terms of regularity on the initial data.
\end{remarka}
\section{Proof of theorem \ref{theo1}}
\label{section2}
We wish to use the proposition \ref{Solo}, that is why we are going to consider an approximation of the system (\ref{1}). To do this we look at viscosity coefficients of the form $\mu_n(\rho)=\max(\frac{1}{n},\mu(\rho))$ with $n\in\mathbb{N}^*$, then the viscosity coefficient $\mu_n$ verifies the assumption (\ref{2.6}) and we can apply the proposition \ref{Solo}. In addition we smooth out the initial data and we work with initial data of the form:
 $$(\rho_0^n-\bar{\rho},u_0^n)=((\rho_0-\bar{\rho})*K_n ,u_0*K_n),$$ with $K_n$ a regularizing kernel. Using the proposition \ref{Solo} there exist approximated strong solutions $(\rho_n,u_n)$ on a maximum time interval $(0,T_n)$. Furthermore since $(\rho_0^n-\bar{\rho},u_0^n)=((\rho_0-\bar{\rho})*K_n ,u_0*K_n)$ is in any $H^s(\R)$ with $1\leq s<+\infty$, therefore we know that the solution $(\rho_n,u_n,v_n)$ is regular and verify in particular $v_n\in L^\infty((0,T_n),L^\infty(\R))$.\\
We are interested now in proving that $T_n=+\infty$ for any $n\in\mathbb{N}$. It will be sufficient if we show that we have for any $T\in[0,T_n]$ (if $T_n<+\infty$) and any $n\in\mathbb{N}$:
\begin{equation}
\begin{aligned}
&\|\rho_n(T,\cdot)-\bar{\rho}\|_{L^\infty(0,T;H^1(\R))}\leq C(T),\\
&\|u_n\|_{L^\infty(0,T;H^1(\R))\cap L^2(0,T,H^2(\R))}+\|\p_t u_n\|_{L^2(0,T,H^2(\R))}\leq C(T),\\
&0<\alpha (T)\leq\rho_n(T,x)\leq\beta(T)<+\infty\;\;\forall x\in\times\R,
\end{aligned}
\label{2.556}
\end{equation}
with $\alpha,\beta$ and $C$ continuous.\\
In addition it implies that for $n\geq \frac{1}{\alpha(T)^\alpha}$, $(\rho_n,u_n)$ is a solution of the system (\ref{1}) on $[0,T]$ for initial data $(\rho_0^n-\bar{\rho},u_0^n)=((\rho_0-\bar{\rho})*K_n ,u_0*K_n)$. Furthermore via the proposition \ref{Solo} we know that this solution is unique.
\\
Now by standard compactness argument, we can prove that the sequence $(\rho_n,u_n)_{n\in\mathbb{N}^*}$ converges up to a subsequence in the sense of the distribution to a global solution $(\rho_1,u_1)$ of (\ref{1}) with initial data $(\rho_0,u_0)$ as in the theorem \ref{theo1} and verifying (\ref{2.55}) and (\ref{2.555}).\\
The uniqueness of the solution will be ensure by the proposition \ref{Solo}. Indeed from  (\ref{2.555}) we have on any interval $t\in[0,T]$ with $T>0$:
$$\mu(\rho_1(t,x))\geq \alpha(T)^{\alpha}\geq c>0\;\;\forall x\in\R,$$
we can then use the uniqueness part of the proposition \ref{Solo} to ensure that the solution $(\rho_1,u_1)$ is unique.\\
Now we are going to focus us on the main difficulty of the proof which consists in proving (\ref{2.556}). In order to simplify the notation in the sequel, we assume that $(\rho,u)$ is a strong solution of (\ref{1}) on the maximal time interval $(0,T_0)$ with initial data $(\rho_0,u_0)$ (it is possible to prove such result if we extend the proposition \ref{Solo} to the case $\mu(\rho)=\rho^{\alpha}$, to do this it suffices to follow the proof of \cite{JDE,JDE2,M3AS} in the case $N\geq 2$). We are going in particular to prove (\ref{2.556}) for $(\rho,u)$ with the maxiimal time interval $(0,T_0)$ and not $(\rho_n,u_n)$ with the maximal time interval $(0,T_n)$ with $n\in\mathbb{N}^*$.\\
In any case the proof does not change when we apply exactly the same argument to $(\rho_n,u_n)$ with $n\in\mathbb{N}^*$. Indeed when we work directly with $(\rho,u)$, we consider always the worth case since $\mu(\rho)$ is degenerate.\\
Let us give an example, typically in order to prove that $v_n$ remains bounded in $L^\infty$ norm, we will need in the  to estimate $\|\frac{P(\rho_n)}{\mu(\rho_n)\rho_n^{\beta}}\|_{L^{\infty}_{t,x}}$ with $\beta=\frac{1}{2}+\e$, $\e>0$. We observe in particular that:
$$\frac{P(\rho_n)}{\mu_n(\rho_n)\rho_n^{\beta}}\leq \frac{P(\rho_n)}{\mu(\rho_n)\rho_n^\beta}.$$
\subsection{Entropy inequalities}
We are going to follow the arguments of \cite{MV}. Let us define the well-known relative entropy, for any functions $U=\biggl(\begin{aligned}
&\rho\\
&\rho u
\end{aligned}\biggl)$ and $\bar{U}=\biggl(\begin{aligned}
&\bar{\rho}\\
&0
\end{aligned}\biggl)$ we set:
$$
\begin{aligned}
{\cal H}(U / \bar{U})&={\cal H}(U)-{\cal H}(\widetilde{U})-D {\cal H}(\bar{U})(U-
\bar{U}),\\
&=\rho (u-\bar{u})^2+p(\rho / \bar{\rho}),
\end{aligned}
$$
where $p(\rho / \bar{\rho})$ is the relative entropy associated to $\frac{1}{\gamma-1}\rho^{\gamma}$:
$$p(\rho / \bar{\rho})=\frac{1}{\gamma-1}\rho^{\gamma}-\frac{1}{\gamma-1}\bar{\rho}^{\gamma}
-\frac{\gamma}{\gamma-1}\bar{\rho}^{\gamma-1}(\rho-\bar{\rho}).$$
Let us mention that since $p$ is convex, the function $p(\rho / \bar{\rho})$ remains positive for every $\rho$ and $p(\rho / \bar{\rho})=0$ if and only if $\rho=\bar{\rho}$.\\
Mellet and Vasseur in \cite{MV} have obtained the following entropy inequalities.
\begin{lem}
Since $(\rho_0, u_0)$ verifies:
\begin{equation}
\int_{\R}{\cal H}(U_0 / \bar{U})dx=\int_{\R}[\rho_0\frac{(u_0)^2}{2}+p(\rho_0 /\bar{\rho})]dx<+\infty,
\label{3.6}
\end{equation}
we have for every $T\in(0,T_0)$, there exists a positive constant $C(T)$ such that:
\begin{equation}
\sup_{[0,T]}\int_{\R}[\rho\frac{(u)^2}{2}+p(\rho/ \bar{\rho})]dx+\int^T_0\int_{\R}\mu(\rho)(\p_x u)^2dx\,dt\leq C(T).
\label{3.7}
\end{equation}
The constant $C(T)$ depends only on $T>0$, $\bar{U}$, the initial value $U_0$, $\gamma$ and on the constant $C$ appearing in (\ref{2.4}).
\label{lemma3.1}
\end{lem}
\begin{remarka}
Let us point out that when both $\bar{\rho}$ and $\rho_0$ are bounded above and below away from zero, we can prove that:
\begin{equation}
p(\rho_0/ \bar{\rho})\leq C(\rho_0-\bar{\rho})^2.
\label{3.6bis}
\end{equation}
In particular (\ref{3.6bis}) is verified under the assumptions of theorem \ref{theo1}.
\end{remarka}
\begin{lem}
For any $T\in(0,T_0)$, there exists $C(T)$ such that the following inequality holds:
\begin{equation}
\begin{aligned}
&\sup_{[0,T]}\int[\frac{1}{2}\rho|u+\p_x(\va(\rho))|^2+p(\rho/ \bar{\rho})|dx\\
&\hspace{4cm}+\int^T_0\int_{\R}\p_x(\va(\rho))\p_x(\rho^\gamma) dx dt\leq C(T),
\end{aligned}
\label{3.12}
\end{equation}
and $\va$ verifying
$\va'(\rho)=\frac{\mu(\rho)}{\rho^2}$. The constant $C(T)$ depends only on $T>0$, $(\bar{\rho},0)$, the initial data $U_0$, $\gamma$ and on the constant $C$ on (\ref{2.4}).
\label{lemma3.2}
\end{lem}
\subsection{Proof of theorem \ref{theo1}}
\subsubsection{A priori estimates}
Since the initial datum $(\rho_0,u_0)$ satisfies (\ref{2.5}), we have:
$$\int\rho_0( u_0)^2 dx<+\infty,\;\;\int p(\rho_0/ \bar{\rho})  dx<+\infty\;\;\mbox{and}\;\;\int \rho_0|\p_x(\va(\rho_0))|^2 dx<+\infty.$$
Using lemmas \ref{lemma3.1} and \ref{lemma3.2} it yields for any $T<T_0$ ($L^\gamma_2$ is the Orlicz space, see \cite{Lio98} in the Appendix A):
\begin{equation}
\begin{aligned}
&\|\sqrt{\rho}u\|_{L^\infty((0,T), L^2)}\leq C(T),\\
&\|\rho\|_{L^\infty((0,T),L^1_{loc}\cap L^\gamma_{loc})}\leq C(T),\\
&\|\rho-\bar{\rho}\|_{L^\infty((0,T),L^\gamma_2)}\leq C(T),\\
&\|\sqrt{\mu(\rho)}\p_x u\|_{L^2((0,T),L^2)}\leq C(T),
\end{aligned}
\label{4.1}
\end{equation}
and:
\begin{equation}
\begin{aligned}
&\|\sqrt{\rho}\p_x\va(\rho)\|_{L^\infty((0,T),L^2)}\leq C(T),\\
&\|\p_x\rho^\gamma\cdot\p_x\va(\rho)\|_{L^2((0,T),L^2)}\leq C(T),
\end{aligned}
\label{4.2}
\end{equation}
with $C$ a continuous function.
\subsubsection{Uniform bounds on the density}
The first proposition shows that the density is bounded by above.
\begin{proposition}
For every $T\in (0,T_0)$ there exists $\beta(\cdot)$ and $C(\cdot)$ continuous functions such that :
$$0\leq \rho(T,x)\leq \beta(T)\;\;\forall x\in\R,$$
and
$$\|\rho-\bar{\rho}\|_{L^\infty((0,T),H^1(\R))}\leq C(T).$$
\label{proposition4.2.2a}
\end{proposition}
{\bf Proof:} Since $\mu(\rho)=\rho^\alpha$ we know via the estimate (\ref{4.2}) that $\p_x\rho^{\alpha-\frac{1}{2}}$ belongs to $L_T^\infty(L^2)$ for any $T\in(0,T_0]$. In addition we know that $\rho-\bar{\rho}$ is in $L_T^\infty(L_\gamma^2)$ using (\ref{4.1}).  For $\e>0$ there exists $C,C_1>0$ such that:
\begin{equation}
C_1|\rho-\bar{\rho}|1_{\{|\rho-\bar{\rho}|\leq\e\}}\leq |\rho^{\alpha-\frac{1}{2}}-\bar{\rho}^{\alpha-\frac{1}{2}}|1_{\{|\rho-\bar{\rho}|\leq\e\}}\leq C|\rho-\bar{\rho}|1_{\{|\rho-\bar{\rho}|\leq\e\}}.
\label{suputi}
\end{equation}
It implies since $\rho-\bar{\rho}$ belongs to $L^\infty(\R^+,L^\gamma_2(\R))$ that there exists $C>0$ such that:
\begin{equation}
\|\big(\rho^{\alpha-\frac{1}{2}}-\bar{\rho}^{\alpha-\frac{1}{2}})1_{\{|\rho-\bar{\rho}|\leq\e\}}\|_{L^2(\R)}\leq C\|\rho-\bar{\rho}\|_{L^\gamma_2(\R)}.
\label{aestim1}
\end{equation}
We have now:
$$(\rho^{\alpha-\frac{1}{2}}-\bar{\rho}^{\alpha-\frac{1}{2}})^2=\rho^{\alpha-\frac{1}{2}}(\rho^{\alpha-\frac{1}{2}}-\bar{\rho}^{\alpha-\frac{1}{2}})+\bar{\rho}^{\alpha-\frac{1}{2}}(\bar{\rho}^{\alpha-\frac{1}{2}}-\rho^{\alpha-\frac{1}{2}}).$$ 
We deduce since the measure of ${\{|\rho-\bar{\rho}|\geq\e\}}$ is finite ( it suffices to use the Tchebytchev inequality and the fact that $\rho-\bar{\rho}$ belongs to $L^{\gamma}_{2}(\R)$):
\begin{equation}
\begin{aligned}
&\|\big(\rho^{\alpha-\frac{1}{2}}-\bar{\rho}^{\alpha-\frac{1}{2}}\big)1_{\{|\rho-\bar{\rho}|>\e\}}   \|^2_{L^2(\R)}\leq\big(\| \rho\|_{L^\infty(\R)}^{\alpha-\frac{1}{2}}+\|\bar{\rho}^{\alpha-\frac{1}{2}}\|_{L^\infty}\big) \|(\rho^{\alpha-\frac{1}{2}}-\bar{\rho}^{\alpha-\frac{1}{2}})1_{\{|\rho-\bar{\rho}|>\e\}} \|_{L^1(\R)}\\
&\leq\big(\| \rho\|_{L^\infty(\R)}^{\alpha-\frac{1}{2}}+\|\bar{\rho}^{\alpha-\frac{1}{2}}\|_{L^\infty}\big) \|(\rho^{\alpha-\frac{1}{2}}-\bar{\rho}^{\alpha-\frac{1}{2}})1_{\{|\rho-\bar{\rho}|>\e\}} \|_{L^2 (\R)}|\{{|\rho-\bar{\rho}|\geq\e\}}|^{\frac{1}{2}}\\
&\leq\big(\| \rho\|_{L^\infty(\R)}^{\alpha-\frac{1}{2}}+\|\bar{\rho}^{\alpha-\frac{1}{2}}\|_{L^\infty}\big) \|(\rho^{\alpha-\frac{1}{2}}-\bar{\rho}^{\alpha-\frac{1}{2}})1_{\{|\rho-\bar{\rho}|>\e\}} \|_{L^2(\R)}\frac{\|\rho-\bar{\rho}\|_{L^2_\gamma} ^{\frac{\gamma}{2}}}{\e^{\frac{\gamma}{2}}}.
\end{aligned}
\label{aestim2}
\end{equation}
The fact that $\alpha-\frac{1}{2}>0$ plays a crucial role here since we can bound 
 $\rho^{\alpha-\frac{1}{2}}$ by $\|\rho\|_{L^\infty}^{\alpha-\frac{1}{2}}$. Now by Sobolev embedding, Gagliardo-Nirenberg inequality , (\ref{aestim1}) and (\ref{aestim2})
, we deduce that there exists $C,C'>0$ sufficiently large such that:
\begin{equation}
\begin{aligned}
&\|\rho^{\alpha-\frac{1}{2}}-\bar{\rho}^{\alpha-\frac{1}{2}}  \|^2_{H^1(\R)}\leq\big(\|\p_x \rho^{\alpha-\frac{1}{2}}\|_{L^2(\R)}+\|\p_x \bar{\rho}^{\alpha-\frac{1}{2}}\|_{L^2(\R)}+\|(\rho^{\alpha-\frac{1}{2}}-\bar{\rho}^{\alpha-\frac{1}{2}})1_{\{|\rho-\bar{\rho}|\leq\e\}}\|_{L^2(\R)}\\
&+\|\big(\rho^{\alpha-\frac{1}{2}}-\bar{\rho}^{\alpha-\frac{1}{2}}\big)1_{\{|\rho-\bar{\rho}|\geq\e\}}   \|_{L^2(\R)}\big)^2\\[2mm]
&\leq C\big(\|\p_x \rho^{\alpha-\frac{1}{2}}\|^2_{L^2(\R)}+\|\p_x \bar{\rho}^{\alpha-\frac{1}{2}}\|^2_{L^2(\R)}+C\|\rho-\bar{\rho}\|^2_{L^\gamma_2(\R)}\\
&\hspace{4cm}+\big(\| \rho^{\alpha-\frac{1}{2}}\|_{L^\infty(\R)}+\|\bar{\rho}^{\alpha-\frac{1}{2}}\|_{L^\infty}\big) \|(\rho^{\alpha-\frac{1}{2}}-\bar{\rho}^{\alpha-\frac{1}{2}}) \|_{H^1(\R)}\frac{\|\rho-\bar{\rho}\|_{L^2_\gamma} ^{\frac{\gamma}{2}}}{\e^{\frac{\gamma}{2}}} \big)\\[2mm]
&\leq C\big(\|\p_x \rho^{\alpha-\frac{1}{2}}\|^2_{L^2(\R)}+\|\p_x \bar{\rho}^{\alpha-\frac{1}{2}}\|^2_{L^2(\R)}+C\|\rho-\bar{\rho}\|^2_{L^\gamma_2(\R)}\\
&+\biggl((\|\p_x( \rho^{\alpha-\frac{1}{2}}-\bar{\rho}^{\alpha-\frac{1}{2}})\|^{\frac{1}{2}}_{L^2(\R)}\| \rho^{\alpha-\frac{1}{2}}-\bar{\rho}^{\alpha-\frac{1}{2}}\|^{\frac{1}{2}}_{L^2(\R)}+2\|\bar{\rho}^{\alpha-\frac{1}{2}}\|_{L^\infty}\biggl)\\
&\hspace{7cm}\times \|(\rho^{\alpha-\frac{1}{2}}-\bar{\rho}^{\alpha-\frac{1}{2}}) \|_{H^1(\R)}\frac{\|\rho-\bar{\rho}\|_{L^2_\gamma} ^{\frac{\gamma}{2}}}{\e^{\frac{\gamma}{2}}} \big)
\end{aligned}
\label{estim3}
\end{equation}
Since by energy estimate $\p_x\rho^{\alpha-\frac{1}{2}}$ belongs to $L^\infty([0,T_0],L^2)$ and $\p_x  \bar{\rho}^{\alpha-\frac{1}{2}}$ is in $L^2$, from bootstrap argument we deduce that $\rho^{\alpha-\frac{1}{2}}-\bar{\rho}^{\alpha-\frac{1}{2}}$ belongs to $L^\infty_T(H^1(\R))$ for any $T\in (0,T_0)$ with:
\begin{equation}
\|\rho^{\alpha-\frac{1}{2}}(T,\cdot)-\bar{\rho}^{\alpha-\frac{1}{2}} (T,\cdot)\|_{H^1(\R))}\leq C(T),
\label{H1}
\end{equation}
with $C$ continuous.
By Sobolev embedding, using the fact that $\alpha-\frac{1}{2}>0$ and since $\bar{\rho}$ belongs to $L^\infty$ we deduce that $\rho$ is bounded in $L^\infty_T(L^\infty(\R))$ for any $T\in(0,T_0)$:
\begin{equation}
\|\rho(T,\cdot)\|_{L^\infty}\leq\beta(T),
\label{normLinf}
\end{equation}
with $\beta$ continuous.\\
Now combining the fact that $\p_x\rho^{\alpha-\frac{1}{2}}$ belongs to $L^\infty([0,T_0], L^2)$ and the fact that $\rho$ is bounded in $L^\infty([0,T_0], L^\infty)$ we deduce that $\p_x\rho$ is bounded in  $L^\infty([0,T_0], L^2)$ (because $\alpha-\frac{1}{2}\leq 1$). Let us prove now that $\rho-\bar{\rho}$ is bounded in $L^\infty([0,T_0], L^2)$, from (\ref{suputi}) we get for $C_1>0$:
\begin{equation}
\begin{aligned}
&\|\rho-\bar{\rho}\|_{L^2}^2\leq \frac{1}{C_1^2}\| \rho^{\alpha-\frac{1}{2}}-\bar{\rho}^{\alpha-\frac{1}{2}}\|_{L^2}^2+ |{\{|\rho-\bar{\rho}|\geq\e\}}|\,\|\rho-\bar{\rho}\|_{L^\infty}^2\\
&\leq \frac{1}{C_1^2}\| \rho^{\alpha-\frac{1}{2}}-\bar{\rho}^{\alpha-\frac{1}{2}}\|_{L^2}^2+ \frac{\|\rho-\bar{\rho}\|_{L^\gamma_2}^\gamma}{\e^\gamma}\,\|\rho-\bar{\rho}\|_{L^\infty}^2.
\end{aligned}
\end{equation}
From (\ref{H1}), (\ref{normLinf}) and energy estimate we deduce that there exists a continuous function $C$ such that for any $T\in (0,T_0)$:
\begin{equation}
\|\rho(T,\cdot)-\bar{\rho} (\cdot)\|_{H^1(\R))}\leq C(T).
\label{H2}
\end{equation}
$\blacksquare$
\\
\\
We are going now to prove that $\rho^{\beta}u$ is bounded in $L^2_T(L^\infty)$ for $\beta$ suitably chosen.
\begin{proposition}
Let $\e>0$ small enough and $\alpha\in (\frac{1}{2},1]$, we have for $\beta=\frac{1}{2}+\e$:
\begin{equation}
\rho^{\beta} u\in L^2([0,T_0],L^\infty).
\label{impimp}
\end{equation}
\label{proposition4.2.2b}
\end{proposition}
{\bf Proof:} Multiplying the momentum equation of (\ref{1}) by $u|u|^p$ with $p>0$ 
as in \cite{MV06} and integrating by parts on $[0,T_0)\times\R$ we can prove that $\rho^{\frac{1}{2+p}}u$ belongs to $L^\infty([0,T_0],(L^{2+p}(\R))$ for any $0<p<+\infty$. Indeed we have for any $T\in (0,T_0)$:
\begin{equation}
\begin{aligned}
&\frac{1}{p+2}\int_{\R}\rho(T,x)|u|^{2+p}(T,x) dx+(p+1)\int^T_0\int_{\R}\mu(\rho(t,x))(\p_x u(t,x))^2|u|^{p} (t,x)dx dt\\
&+\int^T_0\int_{\R}\p_x P(\rho)(t,x) u(t,x)|u|^{p}(t,x) dt dx=\frac{1}{p+2}\int_{\R}\rho_0(x)|u_0|^{2+p}(x) dx.
\end{aligned}
\label{ihyp1}
\end{equation}
We recall that since $u_0$ belongs to $L^\infty(\R)\cap L^2(\R)$ then $u_0$ is any $L^{p+2}(\R)$. Next we have by integration by parts and Young inequality:
\begin{equation}
\begin{aligned}
&|\int^T_0\int_{\R}\p_x P(\rho)(t,x) u(t,x)|u|^{p}(t,x) dt dx|=|(p+1)\int^T_0\int_{\R}\rho^\gamma(t,x)|u|^{p}(t,x)\p_x u(t,x) dt dx|\\
&\leq \frac{p+1}{2}\int^T_0\int_{\R}\rho^{\alpha}(t,x)(\p_x u(t,x))^2 |u|^{p}(t,x) dt dx+\frac{p+1}{2}\int^T_0\int_{\R} \rho^{2\gamma-\alpha}(t,x) |u|^p (t,x) dt dx .
\end{aligned}
\label{ihyp2}
\end{equation}
Next from H\"older inequality and by interpolation, we get for $p\geq 2$ and since $2\gamma-\alpha-1\geq 0$:
\begin{equation}
\begin{aligned}
&\int^T_0\int_{\R} \rho^{2\gamma-\alpha}(t,x) |u|^p (t,x) dt dx\leq \|\rho\|_{L^\infty([0,T]\times\R)}^{2\gamma-\alpha-1} \int^T_0\int_{\R} \rho(t,x) |u|^p (t,x) dt dx \\
&\leq  \|\rho\|_{L^\infty([0,T]\times\R)}^{2\gamma-\alpha-1} \int^T_0\|\rho^{\frac{1}{p+2}}u(t,\cdot)\|_{L^{p+2}}^{\frac{(p-2)(p+2)}{p}}\|\sqrt{\rho}u(t,\cdot)\|_{L^2}^{\frac{4}{p}}dt\\
&\leq  \|\rho\|_{L^\infty([0,T]\times\R)}^{2\gamma-\alpha-1} \int^T_0(1+\|\rho^{\frac{1}{p+2}}u(t,\cdot)\|_{L^{p+2}}^{p+2})\|\sqrt{\rho}u(t,\cdot)\|_{L^2}^{\frac{4}{p}}dt
\end{aligned}
\label{ihyp3}
\end{equation}
Plugging (\ref{ihyp3}) and (\ref{ihyp2}) in (\ref{ihyp1}) and using Gronwall lemma we deduce that $\rho^{\frac{1}{p}}u$ belongs to $L^\infty([0,T_0], L^p)$ for any $p\geq 4$ and by interpolation for any $p\geq 2$.\\
We recall now that since:
\begin{equation}
\p_x (\rho^\beta u)=\rho^{\beta-\frac{\alpha}{2}} \rho^{\frac{\alpha}{2}}\p_x u+\frac{\beta}{\alpha-\frac{1}{2}}\rho^{\beta-\alpha+\frac{1}{2}}u\p_x(\rho^{\alpha-\frac{1}{2}}),
\label{usetech}
\end{equation}
taking $\beta=\frac{1}{2}+\e$ with $0<\e<\frac{1}{4}$, we have:
$$\p_x (\rho^{\beta} u)= \rho^{\frac{1-\alpha}{2}+\e}\rho^{\frac{\alpha}{2}}\p_x u+\frac{\beta}{\alpha-\frac{1}{2}}\rho^{1+\e-\alpha}u\,\p_x(\rho^{\alpha-\frac{1}{2}}).$$
We note now that $\rho^{1+\e-\alpha}u$ belongs to $L^\infty([0,T_0],L^{\frac{1}{\e}}(\R))$ because $\rho$ belongs to $L^\infty([0,T_0],L^\infty)$ and $\rho^{\frac{1}{p}}u$ is in $L^\infty([0,T_0],L^p)$ for any $p\geq 2$. It implies via the estimate (\ref{4.2}), (\ref{normLinf}) and since $\alpha\leq 1$ that $\p_x(\rho^\beta u)$ belongs to $L^2_T(L^2(\R))+L^\infty_T(L^p(\R))$ (for any $T\in(0,T_0)$) which is embedded in $L^2_T(L^2(\R)+L^p(\R))$ with $\frac{1}{p}=\frac{1}{2}+\e$.
By the Riesz Thorin theorem it implies that the Fourier transform ${\cal F}(\p_x(\rho^\beta u))$ is in $L^2_T(L^2(\R)+L^{p'}(\R))$ with $\frac{1}{p}+\frac{1}{p'}=1$. In particular we deduce from H\"older inequality that ${\cal F}(\rho^\beta u)1_{\{|\xi|\geq 1\}}$ is in $L^2_T(L^1(\R))$ for any $T\in [0,T_0]$. Since $\rho^\beta u=\rho^{\e}\sqrt{\rho}u$, we obtain from H\"older inequality that $\rho^\beta u$ is in $L^2_T(L^2)$ for any $T\in [0,T_0]$. From Plancherel theorem we can prove that ${\cal F}(\rho^\beta u)1_{\{|\xi|\leq 1\}}$ is in $L^2_T(L^1(\R))$. It gives that  ${\cal F}(\rho^\beta u)$ is in $L^2_T(L^1(\R))$. We thus get that $\rho^\beta u$ belongs to  $L^2([0,T_0],L^\infty(\R))$. $\blacksquare$\\
\\
The following proposition is the most crucial point of the proof. We show that the density is bounded by below.
\begin{proposition}
For every $T>0$ there exists a continuous function $\alpha$ and $c>0$ such that for all $T<T_0$:
$$0<c\leq\alpha(T)\leq \rho(T,x)\;\;\forall x\in\R.$$
\label{proposition4.2}
\end{proposition}
{\bf Proof:} Since $(\rho, u)$ is a regular solution \footnote{If necessary we can use a regularizing process by looking at initial data of the form $(\rho_0^n-\bar{\rho},u_0^n)=((\rho_0-\bar{\rho})*K_n ,u_0*K_n)$ with $K_n$ a regularizing kernel and passing after to the limit by compactness.} which verifies the system (\ref{1}) on $(0,T_0)$, it is also solution of the system (\ref{11}) on $(0,T_0)$ with:
\begin{equation}
\begin{cases}
\begin{aligned}
&\p_t\rho-\p_{x}(\frac{\mu(\rho)}{\rho}\p_x\rho)+\p_x(\rho v)=0,\\
&\rho\p_t v+\rho u\p_x v+\p_x P(\rho)=0.
\end{aligned}
\end{cases}
\label{moment54}
\end{equation}
We are interested in proving that $v$ remains bounded in $L^\infty ([0,T_0],L^\infty)$. 
Let us multiply the momentum equation of (\ref{moment54}) by $v|v|^p$ (with $p\geq 0$) and integrate over $(0,T)\times\R$ (with $T\in(0,T_0)$), we get using the fact that $\p_x P(\rho)=\gamma \frac{\rho^{\gamma+1}}{\mu(\rho)}(v-u)$:
\begin{equation}
\begin{aligned}
&\frac{1}{p+2}\int_{\R}\rho(T,x)|v(T,x)|^{p+2} dx+\gamma \int^T_0\int_{\R}\frac{\rho^{\gamma+1}(s,x)}{\mu(\rho(s,x))}|v(s,x)|^{p+2} dx\,ds\\
&=\frac{1}{p+2}\int_{\R}\rho_0(x)|v_0(x)|^{p+2} dx+\gamma \int^T_0\int_{\R}\frac{\rho^{\gamma+1}(s,x)}{\mu(\rho(s,x))}u(s,x)v(s,x)|v(s,x)|^{p} dx\,ds.
\end{aligned}
\label{uti1}
\end{equation}
From the Theorem \ref{theo1}, it is easy to observe that $\rho_0^{\frac{1}{q}}v_0$ belongs to any $L^q$ with $q\geq 2$. From H\"older inequality we deduce that:
\begin{equation}
\begin{aligned}
&|\int^T_0\int_{\R}\frac{\rho^{\gamma+1}(s,x)}{\mu(\rho(s,x))}u(s,x)v(s,x)|v(s,x)|^{p} dx\,ds|\\
&=|\int^T_0\int_{\R}\frac{\rho^{\gamma+\frac{1}{p+2}}(s,x)}{\mu(\rho(s,x))}u(s,x)\rho^{\frac{p+1}{p+2}}(s,x) v(s,x)|v(s,x)|^{p} dx\,ds\\
&\leq \int^T_0\|\rho^{\frac{1}{p+2}}v(s,\cdot)\|_{L^{p+2}}^{p+1}\|\frac{\rho^{\gamma+\frac{1}{p+2}}(s,\cdot)}{\mu(\rho(s,\cdot))}u(s,\cdot)\|_{L^{p+2}} ds.
\end{aligned}
\label{uti2}
\end{equation}
It suffices now to control the norm $\|\frac{\rho^{\gamma+\frac{1}{p+2}}(s,\cdot)}{\mu(\rho(s,\cdot))}u(s,\cdot)\|_{L^{p+2}}$. Since $\mu(\rho)=\rho^\alpha$ we have to estimate $\|\rho^{\gamma-\alpha +\frac{1}{p+2}}(s,\cdot) u(s,\cdot)\|_{L^{p+2}}$ and using H\"older inequality we have with $\beta=\frac{1}{2}+\e$ for $\e>0$ small enough:
\begin{equation}
\begin{aligned}
&\|\rho^{\gamma-\alpha -\beta+\frac{1}{p+2}}(s,\cdot) \rho^\beta (s,\cdot)u(s,\cdot)\|^{p+2}_{L^{p+2}}\\
&\hspace{3cm}\leq \|\rho^\beta u(s,\cdot)\|_{L^\infty}^p \|\sqrt{\rho}u(s,\cdot)\|_{L^2}^2\|\rho(s,\cdot)\|_{L^\infty}^{(p+2)(\gamma-\alpha-\frac{p\beta}{p+2})}.
\end{aligned}
\label{tech1}
\end{equation}
We need at this level that $\gamma-\alpha-\frac{p\beta}{p+2}\geq 0$ for any $p\geq 2$ which is true since $\gamma-\alpha-\beta\geq 0$ (indeed in other case we should control the norm $\|\frac{1}{\rho}(s,\cdot)\|_{L^\infty}$ that we precisely wish to estimate). The previous inequality is true since in the Theorem \ref{theo1} we assume that $\gamma\geq \alpha+\frac{1}{2}+\e$.
Plugging (\ref{tech1}) in 
(\ref{uti1}) and using (\ref{uti2}) we have:
\begin{equation}
\begin{aligned}
&\int_{\R}\rho(T,x)|v(T,x)|^{p+2} dx+\gamma (p+2) \int^T_0\int_{\R}\frac{\rho^{\gamma+1}(s,x)}{\mu(\rho(s,x))}|v(s,x)|^{p+2} dx\,ds\\
&\leq \int_{\R}\rho_0(x)|v_0(x)|^{p+2} dx+\gamma(p+2) \int^T_0\big(\|\rho^{\frac{1}{p+2}}v(s,\cdot)\|_{L^{p+2}}^{p+1}\\
&\hspace{3cm}\times\|\rho^\beta u(s,\cdot)\|_{L^\infty}^{\frac{p}{p+2}} \|\sqrt{\rho}u(s,\cdot)\|_{L^2}^{\frac{2}{p+2}}\|\rho(s,\cdot)\|_{L^\infty}^{(\gamma-\alpha-\frac{p\beta}{p+2})}\big)ds.
\end{aligned}
\label{uti3}
\end{equation}
From Young inequality we have:
\begin{equation}
\begin{aligned}
&\int_{\R}\rho(T,x)|v(T,x)|^{p+2} dx+\gamma(p+2) \int^T_0\int_{\R}\frac{\rho^{\gamma+1}(s,x)}{\mu(\rho(s,x))}|v(s,x)|^{p+2} dx\,ds\\
&\leq \int_{\R}\rho_0(x)|v_0(x)|^{p+2} dx+\gamma (p+2) \int^T_0\biggl((1+\|\rho^{\frac{1}{p+2}}v(s,\cdot)\|_{L^{p+2}}^{p+2})\\
&\times\|\rho^\beta u(s,\cdot)\|_{L^\infty}^{\frac{p}{p+2}} \|\sqrt{\rho}u(s,\cdot)\|_{L^2}^{\frac{2}{p+2}}\|\rho(s,\cdot)\|_{L^\infty}^{(\gamma-\alpha-\frac{p\beta}{p+2})}\biggl)ds.
\end{aligned}
\label{uti4}
\end{equation}
Using the Gronwall lemma, we deduce that for any $t\in(0,T_0)$ we have:
\begin{equation}
\begin{aligned}
&\|\rho(t,\cdot)^{\frac{1}{p+2}}v(t,\cdot)\|^{p+2}_{L^{p+2}} \leq \\
&\big( \|\rho_0^{\frac{1}{p+2}}v_0\|^{p+2}_{L^{p+2}}+\gamma (p+2) \int^{T_0}_0\|\rho^\beta u(s,\cdot)\|_{L^\infty}^{\frac{p}{p+2}} \|\sqrt{\rho}u(s,\cdot)\|_{L^2}^{\frac{2}{p+2}}\|\rho(s,\cdot)\|_{L^\infty}^{(\gamma-\alpha-\frac{p\beta}{p+2})}  \,ds\big)\\
&\hspace{0,5cm}\times \exp{\biggl(\gamma(p+2)\int^{T_0}_0\big(\|\rho^\beta u(s,\cdot)\|_{L^\infty}^{\frac{p}{p+2}} \|\sqrt{\rho}u(s,\cdot)\|_{L^2}^{\frac{2}{p+2}}\|\rho(s,\cdot)\|_{L^\infty}^{(\gamma-\alpha-\frac{p\beta}{p+2}))}  \big)\,ds\biggl)}.
\end{aligned}
\label{uti5}
\end{equation}
We deduce that we have for any $p\geq 0$:
\begin{equation}
\begin{aligned}
&\|\rho(t,\cdot)^{\frac{1}{p+2}}v(t,\cdot)\|_{L^{p+2}} \leq \\
&\big( \|\rho_0^{\frac{1}{p+2}}v_0\|^{p+2}_{L^{p+2}}+\gamma (p+2) \int^{T_0}_0\|\rho^\beta u(s,\cdot)\|_{L^\infty}^{\frac{p}{p+2}} \|\sqrt{\rho}u(s,\cdot)\|_{L^2}^{\frac{2}{p+2}}\|\rho(s,\cdot)\|_{L^\infty}^{(\gamma-\alpha-\frac{p\beta}{p+2})}  \,ds\big)^{\frac{1}{p+2}}\\
&\hspace{0,5cm}\times \exp{\biggl(\gamma \int^{T_0}_0\big(\|\rho^\beta u(s,\cdot)\|_{L^\infty}^{\frac{p}{p+2}} \|\sqrt{\rho}u(s,\cdot)\|_{L^2}^{\frac{2}{p+2}}\|\rho(s,\cdot)\|_{L^\infty}^{(\gamma-\alpha-\frac{p\beta}{p+2}))}  \big)\,ds\biggl)}.
\end{aligned}
\label{uti6}
\end{equation}
According to proposition \ref{proposition4.2.2a}, \ref{proposition4.2.2b} and the inequality (\ref{uti6}), we get that for all $p\in [0,+\infty[$:
\begin{equation}
\begin{aligned}
&\|\rho(T,\cdot)^{\frac{1}{p+2}}v(t,\cdot)\|_{L^{p+2}} \leq C(T)\;\;\forall T\in(0,T_0),
\end{aligned}
\label{uti6}
\end{equation}
with $C$ a continuous function independent on $p$. 
\\
At this level, we can prove a $L^\infty_{t,x}$ control on the velocity $v$. Indeed we have the following proposition.
\begin{proposition}\label{prop5}
We have then:
\begin{equation}
	\forall T\in(0,T_0), \;\|v\|_{ L^\infty_T(L^\infty)}\leq C(T),
\label{Lfini}
\end{equation}
with $C$ a continuous function.
\end{proposition}
{\bf Proof:} We recall that for all $t\in (0,T_0)$ the solution $(\rho,v)$ verifies \footnote{This is true if we consider a regularization on the initial data such that $v_0*K_n$ belongs to $H^s$ with $s>\frac{1}{2}$.}:
\begin{equation}
\rho(t,x)\geq  B_t>0\;\;\forall x\in\R\;\;\mbox{and}\;\;\|v(t,\cdot)\|_{L^\infty}\leq C^1_t<+\infty,
\end{equation}
with possibly $B_t\rightarrow_{t\rightarrow T_0} 0$ and $C^1_t\rightarrow_{t\rightarrow T_0} +\infty$.
We observe that $\forall \e>0$ sufficiently small (such that $\e<\frac{\|v(t,\cdot)\|_{L^\infty}}{2}$), 
we have for any $p\geq 2$ and $t\in(0,T_0)$:
\begin{equation}
\begin{aligned}
\|\rho^{\frac{1}{p}}v(t)\|_{L^p}
&\geq \big(\int_{\{x,\;|v(t,x)|\geq \|v(t,\cdot)\|_{L^\infty}-\e\}}
    \rho(t,x) |v|^p (t,x) dx \big)^{\frac{1}{p}}\\  
&\geq \big(\|v(t,\cdot)\|_{L^\infty}-\e\big) 
    \big(\int_{\{x,\;|v(t,x)|\geq \|v(t,\cdot)\|_{L^\infty}-\e\}}
   \rho(t,x)  dx\big)^{\frac{1}{p}}\\
	&\geq \bigl(\|v(t,\cdot)\|_{L^\infty}-\e\bigr) B_{t}^{\frac{1}{p}}
    \,\big|\{x,\;|v(t,x)|\geq\|v(t,\cdot)\|_{L^\infty}-\e\} \big|^{\frac{1}{p}}.
\end{aligned}
\label{norm}
\end{equation}
Since we have $B_{t}>0$ and 
$0<\big| \{x,\;|v(t,x)|\geq \|v(t,\cdot)\|_{L^\infty}-\e \} \big|<+\infty$, 
we can pass to the limit when $p$ goes to $+\infty$ in \eqref{norm}. 
It implies that for any $\e>0$ small enough, we get using (\ref{uti6}):
	\[
    	\|v(t,\cdot)\|_{L^\infty}-\e\leq C(t)\;\;\forall t\in (0,T_0).
    \]
    It concludes the proof of the proposition \ref{prop5}.
$\blacksquare$
\\
\\
In particular we have proved that $v$ belongs to $L^\infty([0,T_0],L^\infty(\R))$.  We can rewrite the mass equation of (\ref{11}) as follows:
\begin{equation}
 \p_t(\frac{1}{\rho})-\p_x(\rho^{\alpha-1}\p_x(\frac{1}{\rho}))+2(\p_x(\frac{1}{\rho})))^2\rho^\alpha+2\p_x(\frac{1}{\rho})v-\p_x(\frac{v}{\rho})=0.
 \label{para14}
 \end{equation}
We deduce using the maximum principle on $[0,T_0)$ (see \cite{L}) that $\frac{1}{\rho}$ is in $L^\infty([0,T_0],L^\infty(\R))$. Indeed we have used the fact that $v$ is in $L^\infty[[0,T_0],L^\infty(\R))$, that there exists $c>0$ such that $\rho^{\alpha-1}(t,x)\geq c>0$ for all $(t,x)\in (0,T_0)\times\R$ (because $\alpha\leq 1$ and $\rho$ belongs to $L^\infty([0,T_0],L^\infty(\R))$) and that $(\p_x(\frac{1}{\rho})))^2\rho^\alpha\geq 0$.  It conclude the proof of the proposition \ref{proposition4.2}.\\
An other possibility consists in proving that $(\frac{1}{\rho}-\frac{1}{\bar{\rho}})$ is bounded in any $L^\infty([0,T_0],L^{p+2} (\R))$ with $p\geq 0$. To do this it suffices again to multiply the equation (\ref{para14}) by $(\frac{1}{\rho}-\frac{1}{\bar{\rho}})|\frac{1}{\rho}-\frac{1}{\bar{\rho}}|^p$ with $p\geq 2$ and integrate over $(0,T)\times\R$ for any $T\in(0,T_0)$. Indeed we have:
\begin{equation}
\begin{aligned}
&\frac{1}{p+2}\int_{\R}|\frac{1}{\rho(T,x)}-\frac{1}{\bar{\rho}(x)}|^{p+2}dx\\
&+(p+1) \int_{0}^T\int_{\R}\rho^{\alpha-1}(t,x)\big(\p_x(\frac{1}{\rho(t,x)})\big)^2|\frac{1}{\rho(t,x)}-\frac{1}{\bar{\rho}(x)}|^{p}dt dx\\
&+2\int_{0}^T\int_{\R} \big(\p_x(\frac{1}{\rho(t,x)})\big)^2\rho^{\alpha-1}(t,x)\,|\frac{1}{\rho(t,x)}-\frac{1}{\bar{\rho}(x)}|^{p}dt dx \leq  \frac{1}{p+2}\int_{\R}|\frac{1}{\rho_0(x)}-\frac{1}{\bar{\rho}(x)}|^{p+2}dx\\
&+(p+1)\int_{0}^T\int_{\R}\rho^{\alpha-1}(t,x)|\p_x(\frac{1}{\rho(t,x)})|\,|\p_x(\frac{1}{\bar{\rho}(x)})|\,|\frac{1}{\rho(t,x)}-\frac{1}{\bar{\rho}(x)}|^{p}dt dx\\
&+2\int_{0}^T\int_{\R} \big(\p_x(\frac{1}{\rho(t,x)})\big)^2\,\frac{\rho^{\alpha}(t,x)}{\bar{\rho}(x)}\,|\frac{1}{\rho(t,x)}-\frac{1}{\bar{\rho}(x)}|^{p}dt dx\\
&+2\int_{0}^T\int_{\R} |\p_x(\frac{1}{\rho(t,x)})|\,|v(t,x)|\,|\frac{1}{\rho(t,x)}-\frac{1}{\bar{\rho}(x)}|^{p+1} dt dx\\
&+(p+1)\int_{0}^T\int_{\R}\frac{|v(t,x)|}{\rho(t,x)} \,|\p_x(\frac{1}{\rho(t,x)}-\frac{1}{\bar{\rho}(x)})|\,|\frac{1}{\rho(t,x)}-\frac{1}{\bar{\rho}(x)}|^p dt dx
\end{aligned}
\label{superimp43}
\end{equation}
It suffices now to proceed by bootstrap and using Gronwall lemma. Indeed we have for example using Young inequality and for $\e,C>0$:
$$
\begin{aligned}
&\int_{0}^T\int_{\R}\rho^{\alpha-1}(t,x)|\p_x(\frac{1}{\rho(t,x)})|\,|\p_x(\frac{1}{\bar{\rho}(x)})|\,|\frac{1}{\rho(t,x)}-\frac{1}{\bar{\rho}(x)}|^{p}dt dx\\
&\leq\frac{\e}{2}\int_{0}^T\int_{\R}\big(\p_x(\frac{1}{\rho(t,x)})\big)^2|\frac{1}{\rho(t,x)}-\frac{1}{\bar{\rho}(x)}|^{p}dt dx\\
&\hspace{5cm}+\frac{1}{2\e}\int_{0}^T\int_{\R}\big(\p_x(\frac{1}{\bar{\rho}(x)})\big)^2\rho^{2\alpha-2}(t,x) |\frac{1}{\rho(t,x)}-\frac{1}{\bar{\rho}(x)}|^{p}dt dx\\
&\leq\frac{\e}{2}\int_{0}^T\int_{\R}\big(\p_x(\frac{1}{\rho(t,x)})\big)^2|\frac{1}{\rho(t,x)}-\frac{1}{\bar{\rho}(x)}|^{p}dt dx\\
&+\frac{C}{2\e}\int_{0}^T\int_{\R}\big(\p_x(\frac{1}{\bar{\rho}(x)})\big)^2(1+|\frac{1}{\rho(t,x)}-\frac{1}{\bar{\rho}(x)}|^{2}+|\frac{1}{\bar{\rho}(x)}|^{2})
|\frac{1}{\rho(t,x)}-\frac{1}{\bar{\rho}(x)}|^{p}dt dx\\
\end{aligned}
$$
We conclude for this term by using Gronwall lemma and the fact that $\p_x(\frac{1}{\bar{\rho}})$, $\frac{1}{\bar{\rho}}$ are in $L^\infty$ and $(\rho-\bar{\rho})$ belongs to $L^\infty([0,T_0],L^2\cap L^\infty)$. In a similar way we have for $\e>0$:
$$
\begin{aligned}
&\int_{0}^T\int_{\R} \big(\p_x(\frac{1}{\rho(t,x)})\big)^2\,\frac{\rho^{\alpha}(t,x)}{\bar{\rho}(x)}\,|\frac{1}{\rho(t,x)}-\frac{1}{\bar{\rho}(x)}|^{p}dt dx \leq\\
&\e\int_{0}^T\int_{\R} \big(\p_x(\frac{1}{\rho(t,x)})\big)^2\,\rho^{\alpha-1}(t,x)\,|\frac{1}{\rho(t,x)}-\frac{1}{\bar{\rho}(x)}|^{p}dt dx\\
&\hspace{3cm}+\int_{0}^T\int_{\{\rho\geq \e\bar{\rho}\}}\big(\p_x(\frac{1}{\rho(t,x)})\big)^2 \frac{\rho^{\alpha}(t,x)}{\bar{\rho}(x)}\,|\frac{1}{\rho(t,x)}-\frac{1}{\bar{\rho}(x)}|^{p}dt dx.
\end{aligned}
$$
The second term on the right hand side is bounded because $\frac{1}{\bar{\rho}}$ is in $L^\infty$, $\rho$ belongs to $L^\infty([0,T_0],L^\infty)$ and $\p_x\rho$ is bounded in $L^\infty([0,T_0],L^2)$. The first term on the right hand side can be absorbed by the terms on the left hand side of (\ref{superimp43}). Let us give a last example, using again Young inequality we have for $\e>0$:
$$
\begin{aligned}
&\int_{0}^T\int_{\R} |\p_x(\frac{1}{\rho(t,x)})|\,|v(t,x)|\,|\frac{1}{\rho(t,x)}-\frac{1}{\bar{\rho}(x)}|^{p+1} dt dx\leq\\
&\frac{\e}{2} \int_{0}^T\int_{\R} |\p_x(\frac{1}{\rho(t,x)})|^2 \,|\frac{1}{\rho(t,x)}-\frac{1}{\bar{\rho}(x)}|^{p} dt dx\\
&\hspace{5cm}+\frac{1}{2\e} \int_{0}^T\|v(t,\cdot)\|_{L^\infty}^2( \int_{\R} |\frac{1}{\rho(t,x)}-\frac{1}{\bar{\rho}(x)}|^{p+2}  dx) dt.
\end{aligned}
$$
Now using the fact that $\p_x\rho^{\alpha-\frac{1}{2}}$ belongs to $L^\infty([0,T_0],L^2)$ we prove easily that $(\frac{1}{\rho}-\frac{1}{\bar{\rho}})$ is bounded in $L^\infty([0,T_0]\times\R)$. Indeed it suffices to observe that $\p_x(\frac{1}{\rho})=-\frac{1}{\alpha-\frac{1}{2}}(\frac{1}{\rho^{\alpha+\frac{1}{2}}}-\frac{1}{\bar{\rho}^{\alpha+\frac{1}{2}}})\p_x(\rho^{\alpha-\frac{1}{2}})-\frac{1}{\alpha-\frac{1}{2}}\frac{1}{\bar{\rho}^{\alpha+\frac{1}{2}}}\p_x(\rho^{\alpha-\frac{1}{2}})$. Since $(\frac{1}{\rho^{\alpha+\frac{1}{2}}}-\frac{1}{\bar{\rho}^{\alpha+\frac{1}{2}}})$ is bounded in any $L^\infty([0,T_0],L^{p+2})$ with $p\geq 2$ we deduce that $\p_x(\frac{1}{\rho})$ belongs to $L^\infty([0,T_0],L^2+L^q)$ with $\frac{1}{q}=\frac{1}{2}+\frac{1}{p+2}$. Now we obtain the result by using similar arguments than the proof of the proposition \ref{proposition4.2.2b}.
$\blacksquare$
\subsubsection{Uniform bounds for the velocity}
\begin{proposition}
\label{proposition3.6}
There exists a continuous function $C$ such that $\forall T\in (0,T_0)$:
$$\|u\|_{L^2((0,T),H^2(\R))}\leq C(T),$$
and:
$$\|\p_t u\|_{L^2((0,T),L^2(\R))}\leq C(T).$$
In particular $u\in C((0,T),H^1(\R))$, $\forall T\in (0,T_0)$.
\label{proposition4.7}
\end{proposition}
{\bf Proof:} This proposition follows \cite{MV}. First, we observe that $u$ is bounded in $L^2((0,T),H^1(\R))$ for any $T\in(0,T_0)$. Indeed since we have
 $\rho\geq\alpha_1>0$ using (\ref{3.7}) it implies that $\p_x u$ is bounded in $L^2((0,T)\times\R)$ and $u$ is bounded in $L^\infty((0,T),L^2(\R))$. Therefore $u$ belongs to $L^2([0,T_0],H^1(\R))$.\\
We deduce via the mass equation that $\p_t \rho$ is bounded in $L^2([0,T_0]\times\R)$ using (\ref{4.2}), propositions \ref{proposition4.2.2a}, \ref{proposition4.2} and Sobolev embedding. Since $\rho-\bar{\rho}$ is bounded in $L^\infty([0,T_0],H^1(\R))$, it yields  using interpolation and Sobolev embedding (see \cite{Ama00} for example) that $\rho$ belongs to $C^{s_0}([0,T_0]\times\R)$ for some $s_0$ in $(0,1)$.\\
Next let us rewrite the momentum equation of (\ref{1}) as follows:
\begin{equation}
\p_t u-\p_x (\frac{\mu(\rho)}{\rho} \p_x u)=-\gamma\rho^{\gamma-2}\p_x\rho+\big(\p_x\va(\rho)-u\big)\p_x u.
\label{4.3}
\end{equation}
In order to obtain new estimates on $u$, we are going to control the right hand side of (\ref{4.3}). The first term $\rho^{\gamma-2}\p_x\rho$ is bounded in $L^\infty([0,T_0],L^2(\R))$ (thanks to propositions \ref{proposition4.2.2a}, \ref{proposition4.2} and estimate (\ref{4.2})). For the last part following \cite{MV}, we write (using H\"older inequality, Gagliardo-Nirenberg inequality and energy estimate):
$$
\begin{aligned}
&\|\big(\p_x\va(\rho)-u\big)\p_x u\|_{L^2((0,T),L^{2}(\R))}\\
&\leq \|\p_x\va(\rho)-u\|_{L^\infty((0,T),L^2(\R))}\|\p_x u\|_{L^2((0,T),L^{\infty}(\R))}\\
&\leq \|\p_x\va(\rho)-u\|_{L^\infty((0,T),L^2(\R))}\|\p_x u\|^{\frac{1}{2}}_{L^2((0,T),L^2(\R))}
\|\p_{xx} u\|^{\frac{1}{2}}_{L^2((0,T),L^2(\R))}\\
&\leq C \|\p_{xx}u\|^{\frac{1}{2}}_{L^2((0,T),L^2(\R))}.
\end{aligned}
$$
Using regularity results for parabolic equation of the form (\ref{4.3}) ( see \cite{L}, we can check here that that the diffusion coefficient is in $C^{s_0}((0,T)\times\R)$) gives for any $T\in(0,T_0)$:
$$
\|\p_tu\|_{L^2((0,T),L^2(\R)}+ \| \p_x u\|_{L^2 ((0,T),H^1(\R))}\leq C  \|\p_x u\|^{\frac{1}{2}}_{L^2((0,T),H^1(\R))}+C,$$
with $C$ depending on $\|u_0\|_{H^1}$ and by bootstrap for any $T\in(0,T_0)$:
 $$\|\p_tu\|_{L^2((0,T),L^2(\R))}+ \| u\|_{L^2((0,T),H^2(\R))}\leq C(T).$$
This concludes the proof of the proposition \ref{proposition3.6}. $\blacksquare$\\
\\
Combining all the previous propositions, we have proved the theorem \ref{theo1}.
 \section{Appendix}
 Let us explain why the system (\ref{1}) is equivalent to the system (\ref{11}) via the change of unknown $v=u+\p_x\va(\rho)$.
  \begin{proposition}
We can formally rewrite the system (\ref{1}) as follows with $v=u+\p_x\va(\rho)$ ($\va'(\rho)=\frac{\mu(\rho)}{\rho^2}$):
\begin{equation}
\begin{cases}
\begin{aligned}
&\p_{t}\rho-\p_x(\frac{\mu(\rho)}{\rho}\p_x\rho)+\p_x(\rho v)=0,\\
&\rho \p_{t}v+\rho u \p_x v+\p_x P(\rho)=0,\\
&(\rho,u)_{/t=0}=(\rho_{0},u_{0}).
\end{aligned}
\end{cases}
\label{0.1a}
\end{equation}
\end{proposition}
{\bf Proof:}
As observed in \cite{cras,para,CPAM,CPAM1}, we are interested in rewriting the system (\ref{1}) in terms of the following unknown $v=u+\p_x \va(\rho)$ where $\va'(\rho)=\frac{\mu(\rho)}{\rho^2}$. We have then by using the mass equation of (\ref{1}):
$$\p_t\rho+\p_x(\rho v)-\p_x(\rho \va'(\rho)\p_x\rho)=0.$$
and it gives:
$$\p_t\rho+\p_x(\rho v)-\p_x(\frac{\mu(\rho)}{\rho}\p_x\rho)=0.$$
From the mass equation in (\ref{1}) we obtain:
\begin{equation}
\begin{aligned}
&\rho \p_t [\p_x(\va(\rho))]+\rho\p_x\big(\frac{\mu(\rho)}{\rho^2}\p_x(\rho u)\big)=0.
\end{aligned}
\label{mass1}
\end{equation}
Next we have:
\begin{equation}
\begin{aligned}
\rho\p_x\big(\frac{\mu(\rho)}{\rho^2}\p_x(\rho u)\big)&=\rho\p_x(\frac{\mu(\rho)}{\rho}\p_x u+u\p_x\va(\rho)),\\
&=\mu(\rho)\p_{xx} u+\rho\p_x(\frac{\mu(\rho)}{\rho})\p_x u+\rho\p_x(u\p_x\va(\rho)),\\
&=\mu(\rho)\p_{xx} u+\p_x\mu(\rho)\,\p_x u-\frac{\mu(\rho)}{\rho}\p_x\rho\,\p_x u+\rho\p_x(u\p_x\va(\rho)),\\
&=\p_x(\mu(\rho)\p_x u)+\rho u\p_{xx}\va(\rho).
\end{aligned}
\label{acitech}
\end{equation}
Now we can rewrite the momentum equation of (\ref{1}) as follows:
\begin{equation}
\rho\p_t u+\rho u\p_x u-\p_x(\mu(\rho)\p_x u)+\p_x P(\rho)=0.
\label{moment23}
\end{equation}
We get now adding the equality (\ref{mass1}) to the momentum equation (\ref{moment23}) and using (\ref{acitech}):
$$
\begin{aligned}
&\rho \p_t v+\rho u\p_x v+\p_x P(\rho)=0.
\end{aligned}
$$
It concludes the proof. $\blacksquare$


\begin{thebibliography}{}
\bibitem{Ama00}
H. Amann. Compact embeddings of vector-valued Sobolev and Besov spaces. \textit{Glas. Mat. Ser. III}, 35(55) (1): 161-177, 2000. Dedicated to the memory of Branko Najman.
\bibitem{BDa}
D. Bresch and B. Desjardins. Existence of global weak solutions for
a 2D Viscous shallow water equations and convergence to the
quasi-geostrophic model. \textit{Comm. Math. Phys.}, 238(1-2): 211-223,
2003.
\bibitem{BD}
D. Bresch and B. Desjardins. Existence of global weak solutions to the Navier- Stokes equations for viscous compressible and heat conducting fluids, \textit{Journal de Math\'ematiques Pures et Appliqu\'es}, Volume 87, Issue 1, January 2007, Pages 57-90.
\bibitem{BrDeZa}
D. Bresch, B. Desjardins, and E. Zatorska. Two-velocity hydrodynamics in fluid mechanics: Part {II} existence of
  global $\kappa$-entropy solutions to compressible Navier-Stokes systems with
  degenerate viscosities. \textit{J. Math. Pures Appl.} (9)104, no. 4, 801--836, 2015.
\bibitem{BJ}
D. Bresch and E. Jabin. Global Existence of Weak Solutions for Compresssible Navier--Stokes Equations: Thermodynamically unstable pressure and anisotropic viscous stress tensor. \textit{Preprint} 	arXiv:1507.04629.
\bibitem{CC70}
S. Chapman and T.G Cowling. The mathematical theory of non-uniform gases. An account of the kinetic theory of viscosity, thermal conduction and diffusion in gases. \textit{Third edition, prepared in co-operation with D. Burnett. Cambridge University Press, London, 1970.}
\bibitem{Daf79} 
C. M. Dafermos. The second law of thermodynamics and stability. \textit{Arch. Rational Mech. Anal.}, 70(2): 167-179; 1979.
\bibitem{Fei04}
E. Feireisl. On the motion of a viscous compressible, and heat conducting fluid. \textit{Indiana Univ. Math. J.}, 53(6): 1705-1738, 2004.
\bibitem{FNP01}
E. Feireisl, A. Novotn\' y and H. Petzelt\`ova. On the existence of globally defined weak solutions to the Navier-Stokes equations. \textit{J. Math. Fluid Mech.}, 3(4): 358-392, 2001.
\bibitem{GP}
J.-F. Gerbeau and B. Perthame. Derivation of viscous Saint-Venant system for laminar shallow water; numerical validation. \textit{Discrete Contin. Dyn. Syst. Ser. B}, 1(1):89-102, 2001.
\bibitem{M3AS}
B. Haspot. Cauchy problem for viscous shallow water equations with a term of capillarity, \textit{Mathematical Models and Methods in Applied Sciences}, {\bf 20} (7) (2010), 1049-1087.
\bibitem{JDE}
B. Haspot. Well-posedness in critical spaces for the system of compressible NavierÐStokes in larger spaces, \textit{Journal of Differential Equations},  251, No. 8. (October 2011), pp. 2262-2295.
\bibitem{CPAM1}
B. Haspot. Hyperbolic Problems: Theory, Numerics, Applications Ñ Proceedings of the 14th International Conference on Hyperbolic Problems held in Padova, June 25-29, 2012, p. 667-674, 2014.
\bibitem{CPAM}
B. Haspot. From the highly compressible Navier-Stokes equations to fast diffusion and porous media equations,  existence of global weak solution for the quasi-solutions. \textit{Journal of Mathematical Fluid Mechanics}. 18(2) (2016), 243-291.
\bibitem{JDE2}
 B. Haspot. Global existence of strong solution for shallow water system with large initial data on the irrotational part, \textit{Journal of Differential Equations} Volume 262, Issue 10, (2017), 4931-4978.
\bibitem{para}
B. Haspot. New formulation of the compressible Navier-Stokes equations and parabolicity of the density. \textit{Preprint} 2015 Hal.
\bibitem{cras}
B. Haspot and E. Zatorska, From the highly compressible Navier-Stokes equations to the Porous
Media equation, rate of convergence. To appear in \textit{Discrete and Continuous Dynamical Systems - Series A}(36), (2016) 3107-3123.
\bibitem{Hof87}
D. Hoff. Global existence for 1D compressible, isentropic Navier-Stokes equations with large initial data. \textit{Trans. Amer. Math. Doc.}, 303(1): 169-181; 1987.
\bibitem{Hof95}
D. Hoff. Global solutions of the Navier-Stokes equations for multidimensional compressible flow with discontinuous initial data. \textit{J. Differential Equations}, 120(1): 215-254, 1995.
\bibitem{Hof98}
D. Hoff. Global solutions of the equations of one-dimensional, compressible flow with large data and forces, and with differing end states. \textit{Z. Angew. Math. Phys.}, 49(5): 774-785, 1998.
\bibitem{HS01}
D. Hoff and J. Smoller. Non-formation of vacuum states for compressible Navier-Stokes equations. \textit{Comm. Math. Phys.}, 216(2): 255-276, 2001.
\bibitem{JXZ}
S. Jiang, Z. Xin and P. Zhang.  Global weak solutions to 1D compressible isentropic Navier-Stokes equations with density dependent viscosity. \textit{Methods and applications of analysis}
Vol. 12, No. 3, pp. 239Ð252, September 2005.
\bibitem{Ka}
Ya. I. Kanel. On a model system of equations of one-dimensional gas motion, \textit{Differ-
entsial'nye Uravneniya} 4 (1968), 721-734.
\bibitem{KS77}
A. V. Kazhikhov and V. V.  Shelukhin. Unique global solution with respect to time of initial-boundary value problems for one-dimensional equations of a viscous gas. \textit{Prikl. Mat. Mech.}, 41(2): 282-291, 1977.
\bibitem{L}
O. A. Ladyzhenskaya, V. A. Solonnikov and N. N. Uraltseva. Linear and quasilinear equations of parabolic type, \textit{AMS translations, Providence}, 1968.
\bibitem{LLX}
H-L. Li, J. Li and Z. Xin. Vanishing of Vacuum States and Blow-up Phenomena
of the Compressible Navier-Stokes Equations. \textit{Commun. Math. Phys.} 281, 401Ð444 (2008).
\bibitem{Li}
J. Li and Z. Xin. Global Existence of Weak Solutions to the
Barotropic Compressible Navier-Stokes Flows with Degenerate
Viscosities. \textit{arXiv:1504.06826} [math.AP].
\bibitem{Lio98}
P-L. Lions. Mathematical topics in fluid mechanics. \textit{Vol. 2, Volume 10 of Oxford Lecture Series in Mathematics and its Applications}. The Clarendon Press Oxford University Press, New-York, 1998. Compressible models, Oxford science publications.
\bibitem{LXY98} T-P. Liu, Z. Xin and T. Yang, Vacuum states for compressible flow.Ê\textit{Discrete Contin. Dynam. Systems}, 4(1): 1-32, 1998.
\bibitem{MN79}
A. Matsumura and T. Nishida. The initial value problem for the equations of motion of compressible viscous and heat-conductive fluids. \textit{Proc. Japan Acas. Ser. A Math. Sci.}, 55(9): 337-342,1979.
\bibitem{MV06}
A. Mellet and A. Vasseur. On the barotropic compressible Navier-Stokes equations. 
\textit{Comm. Partial Differential Equations} 32 (2007), no. 1-3, 431--452.
\bibitem{MV}
A. Mellet and A. Vasseur. Existence and Uniqueness of global strong solutions for one-dimensional compressible Navier-Stokes equations. 
\textit{SIAM J. Math. Anal.} 39 (2007/08), no. 4, 1344--1365.
\bibitem{Mu}
P. B. Mucha. Compressible Navier-Stokes system in 1-D. \textit{Math. Methods Appl. Sci.}, 24(9):607Ð622, (2001).
\bibitem{Mu1}
P. B. Mucha. The Cauchy problem for the compressible Navier-Stokes equation in the $L_p$ framework, \textit{Nonlinear Analysis}, 52(4), 1379-1392 (2003).
\bibitem{OMNM02}
M. Okada, S. Matusu-Necasova and T. Makino. Free boundary problem for the equation of one-dimensional motion of compressible gas with density-dependent viscosity. \textit{Ann. Univ. Ferrara Sez.} VII (N.S); 48:1-20, 2002.
\bibitem{Ser86a}
D. Serre. Solutions faibles globales des \'equations de Navier-Stokes pour un fluide compressible. \textit{C. R. Acad. Sci. Paris S\'er. I Math.}, 303(13): 639-642, 1986.
\bibitem{Ser86b}
D. Serre. Sur l'\'equation monodimensionnelle d'un fluide visqueux, compressible et conducteur de chaleur. \textit{C. R. Acad. Sci. Paris S\'er. I Math.}, 303 (14): 703-706, 1986.
\bibitem{She82}
V. V. Shelukin. Motion with a contact discontinuity in a viscous heat conducting gas. \textit{Dinamika Sploshn. Sredy}, (57): 131-152, 1982.
\bibitem{She83}
V. V. Shelukin. Evolution of a contact discontinuity in the barotropic flow of a viscous gas. \textit{Prikl. Mat. Mekh.} 47(5): 870-872, 1983.
\bibitem{She84}
V. V. Shelukin. On the structure of generalized solutions of the one-dimensional equations of a polytropic viscous gas. \textit{Prikl. Mat. Mekh.}, 48(6): 912-920; 1984.
\bibitem{She86}
V. V. Shelukin. Boundary value problems for equations of a barotropic viscous gas with nonnegative initial density. \textit{Dinamika Sploshn. Sredy}, (74): 108-125, 162-163, 1986.
\bibitem{She98}
V. V. Shelukin. A shear flow problem for the compressible Navier-Stokes equations. \textit{Internat. J. Non-linear Mech.}, 33(2): 247-257, 1998.
\bibitem{Sol76}
V. A. Solonnikov. The solvability of the initial-boundary value problem for the equations of motion of a viscous compressible fluid. \textit{Zap. Naucn. Sem. Leningrad. Otdel. Mat. Inst. Steklov.} (LOMI), 56: 128-142, 197, 1976. Investigations on linear operators and theory of functions, VI.
\bibitem{Vai90}
V. A Vaigant. Nonhomogeneous boundary value problems for equations of a viscous heat-conducting gas. \textit{Dinamika Sploshn. Sredy,} (97): 3-21, 212, 1990.
 \bibitem{V}
A. F. Vasseur and C. Yu. Existence of global weak solutions for 3D degenerate compressible Navier-Stokes equations. \textit{Invent. Math}.  206 (2016), no. 3, 935-974.
\bibitem{V1}
A. F. Vasseur and C. Yu. Global weak solutions to compressible quantum Navier-Stokes equations with damping. \textit{SIAM J. Math. Anal.} 48 (2016), no. 2, 1489-1511.
\bibitem{YYZ01}
T. Yang, Z-a Yao and C. Zhu. Compressible Navier-Stokes equations with density-dependent viscosity and vacuum. \textit{Comm. Partial Differential Equations}, 26(5-6): 965-981, 2001.
\bibitem{YZ02}
T. Yang and C. Zhu. Compressible Navier-Stokes equations with degenerate viscosity coefficients and vacuum. \textit{Comm. Math. Phys.}, 230(2): 329-363, 2002.
\end{thebibliography}
\end{document}